\documentstyle[11pt]{article}
\newcommand{\be}{\begin{equation}}
\newcommand{\ee}{\end{equation}}
\newcommand{\bea}{\begin{eqnarray}}
\newcommand{\eea}{\end{eqnarray}}
\newcommand{\beann}{\begin{eqnarray*}}
\newcommand{\eeann}{\end{eqnarray*}}

\newcommand{\norm}[1]{\mbox{$\|#1\|$}}

\begin{document}

\title{COMPACT PERTURBATIONS OF 
FREDHOLM n-TUPLES, II}  
\author{R{\u{a}}zvan
Gelca}
\maketitle

\hspace{15 mm}  We study compact perturbations 
of Fredholm n-tuples of index zero.
We prove that if the operators in such an n-tuple acting on a 
Hilbert space satisfy  certain functional relations then the n-tuple
cannot be perturbed with compact operators to an invertible one.

\vspace{15 mm}
 
\hspace{15 mm} The study of commuting
 n-tuples of operators has been initiated by
J. L. Taylor in [5] and [6]. Since then several properties of a single 
operator have been generalized to n-tuples. In [2] R. Curto asked if
the fact that a Fredholm operator of index zero acting on a Hilbert
space can be made invertible by adding a compact operator
remains true for commuting pairs. In [3] it has been shown that pairs of the 
form $(T,T)$ with $T$ Fredholm and $indT\neq 0$ cannot be perturbed to 
invertible pairs. The aim of this paper is to extend this result to n-tuples.

\hspace{15 mm} We shall start by reviewing some important facts about commuting
n-tuples. We consider only operators on a certain infinite
dimensional Hilbert space $H$. Following [5] we attach to each commuting
n-tuple $T=(T_1,T_2,\cdots ,T_n)$ a complex of Hilbert spaces
$(K^p(T,H),\delta _T)$, called the Koszul complex, by defining
$K^p(T,H):=H\bigotimes \Lambda ^p$, and $\delta _T:K^p(T,H)\rightarrow
K^{p+1}(T,H)$, $\delta _T:=T_1\otimes E_1+ \cdots +T_n\otimes E_n$,
where $\Lambda ^p=\Lambda ^p[e_1,e_2,\cdots ,e_n]$ are the p-forms
on ${\bf C}^n$ and \( E_i\omega =e_i\omega$, $\omega \in \Lambda ^p,
i=1,2,\cdots ,n$.

\hspace{15 mm} The n-tuple $T$ is called invertible if its
Koszul complex is exact. The spectrum of $T$, denoted
by $\sigma (T)$, is the set of all
\( z=(z_1,z_2,\cdots ,z_n)\in {\bf C}^n\) such that 
\( z-T=(z_1-T_1,z_2-T_2,\cdots ,z_n-T_n)\) is not invertible.
In [5] it is proved that the spectrum is a compact nonvoid
set. For any holomorphic map $f$ on an open neighborhood
of $\sigma (T)$ one can define $f(T)$ (cf. [6]). 
The spectral mapping theorem asserts that \( f(\sigma (T)=
\sigma (f(T))\) (cf. [6]). A general overview of
the properties of the spectrum can be found in [1].

\hspace{15 mm} Let $H^p(T)$, $0\leq p\leq n$, be the cohomology
spaces of the Koszul complex. The n-tuple $T$ is called Fredholm
if all these spaces are finite dimensional, in this case we define
the index of $T$ to be \( indT=\sum_{p=0}^{n}(-1)^pdimH^p(T)$.
It is known that the index is preserved under compact perturbations
(cf. [4]).
Given two commuting tuples $T=(T_1,T_2,\cdots ,T_n)$ and
$T'=(T_1,T_2,\cdots ,T_n,S)$ one has the long exact sequence in
cohomology 

\hspace{10 mm} \( 0\rightarrow H^0(T')\rightarrow H^0(T)
\stackrel{\hat S_0}{\rightarrow} H^0(T)\rightarrow
H^1(T')\rightarrow H^1(T)\rightarrow \cdots \)

\begin{eqnarray*}
\cdots \rightarrow H^{p-1}(T)\rightarrow H^p(T')\rightarrow
H^p(T)\stackrel{\hat{S_p}}{\rightarrow}
H^p(T)\rightarrow \cdots
\end{eqnarray*}  

\noindent where ${\hat S}_p$ is the operator induced by
\( S\otimes 1:K^p(T,H)\rightarrow K^p(T,H)\), $0\leq p\leq n$.
We observe that the tuple $T'$ is invertible if and only if all the
operators ${\hat S}_p$ are isomorphisms. As a consequence of 
this long exact sequence if $T$ is Fredholm then $T'$ is
Fredholm of index zero; this is a method of obtaining
Fredholm n-tuples of index zero. Given a subspace ${\cal H}$
of $H$ we denote by $P_{\cal H}$ the orthogonal projection on
${\cal H}$.
\bigskip

\begin{em}
\hspace{15 mm} {\bf{Lemma.}} Let $T$ be such that for any $n$,
$dim ker T^n<\infty$ and $dimkerT^n \rightarrow \infty$.
If $S$ commutes with $T$ and the sequence $dim(kerS\bigcap 
kerT^n)$, $n\in {\bf{N}}$ is bounded then there exists a 
sequence of nontrivial orthogonal subspaces ${\cal H}_n$
in $H$ such that $P_{{\cal H}_n}S|{\cal H}_n$
 is invertible, and for
every $m$ and $n$, $P_{{\cal H}_n}S|{\cal H}_n$ is
similar to $P_{{\cal H}_m}S|{\cal H}_m$.
\end{em}
\medskip

\hspace{15 mm} {\bf {Proof.}} Let ${\cal K}_n=kerT^n\ominus  kerT^{n-1}$.
Since $dim kerT^n\rightarrow \infty$, the spaces ${\cal K}_n$
are nontrivial. Moreover, the operator

\begin{eqnarray*}
P_{{\cal K}_n}T|{\cal K}_n:{\cal K}_n\rightarrow {\cal K}_{n-1}
\end{eqnarray*}

\noindent  is injective, therefore $dim{\cal K}_n\leq dim{\cal K}_{n-1}$.
This shows that the sequence $dim{\cal K}_n$, $n\in {\bf N}$ is
a decreasing sequence of natural numbers, so it becomes 
stationary. It follows that there exists a number $n_0$ such
that for $n\geq n_0$, the operator $P_{{\cal K}_{n-1}}T|{\cal K}_n$
is an isomorphism. Since for every $n$, $kerT^n\subset kerT^{n+1}$
and the sequence $dim(ker S\bigcap ker T_n)$, $n\in {\bf N}$
is bounded, there exists a number $n_1>n_0$ such that for
$n\geq n_1$, the operator $P_{{\cal K}_n}S|{\cal K}_n$
is injective, hence invertible. Moreover, the operator 
$P_{{\cal K}_n}T|{\cal K}_n$ defines a similarity between 
$P_{{\cal K}_n}S|{\cal K}_n$ and $P_{{\cal K}_{n+1}}S|{\cal K}_{n+1}$
for every $n\geq n_1$. Taking ${\cal H}_n={\cal K}_{n+n_1}$, 
$n\geq 0$, we obtain a sequence of spaces with the desired property.

\newpage

\begin{em}
\hspace{15 mm} {\bf Theorem.} Let $(T_1,T_2,\cdots ,T_n)$ be a commuting
n-tuple with $T_1$ Fredholm and $indT_1$ different from zero. If there
exists for each $k$, $2\leq k\leq n$ an analytic function
of two variables $f_k$ such that

1. $f_k(0,w)=0$ implies $w=0$,

2. $f_k(T_1,T_k)=L_k$, $L_k$ compact,

\noindent  then the n-tuple $(T_1,T_2,\cdots ,T_n)$ cannot be perturbed
with compact operators to an invertible n-tuple.
\end{em}

\vspace{3mm}

\hspace{15 mm} {\bf Proof.} Suppose that such compacts 
$K_1,K_2,\cdots ,K_n$ exist. Denote $S_i=T_i+K_i$.
Then $S_1$ is Fredholm of nonzero index, we may assume
$indS_1>0$. We remark that for every $k$, $2\leq k\leq n$ 
the operator $N_k=f_k (S_1,S_k)$ is compact.
 Consider the analytic function $f:{\bf C}^n \rightarrow {\bf C}^n$,
 $f(z_1,z_2,\cdots ,z_n)$ $=(z_1,f_2(z_1,z_2),\cdots ,f_n(z_1,z_n))$.
Then $f^{-1}(0)=0$, and since $(S_1,S_2,\cdots S_n)$ is invertible,
from the spectral mapping theorem it follows that $(S_1,N_2,\cdots ,N_n)$
is also invertible. Let us show that this is not possible.

\hspace{15 mm} Let us consider 
 $k$ to be the smallest integer with the property that the sequence
$dim(kerS_1^m\bigcap ker N_2 \bigcap \cdots$$  \bigcap kerN_k)$,
$m\in {\bf N}$, is bounded. Such a $k$ exists, for by the
spectral mapping theorem $(S_1^m,N_2,\cdots ,N_n)$ is invertible
for every $m$, hence $kerS_1^m\bigcap kerN_2\bigcap \cdots
 \bigcap kerN_n$ $=0$.
Consider the subspace $H_0=kerN_2\bigcap \cdots \bigcap N_{k-1}$ (in case $k=2$
take $H_0=H$). Since the operators $S_1,N_2,\cdots ,N_k $ commute,
$H_0$ is  invariant for $S_1$ and $N_k$. Moreover, because of the
minimality of $k$, the operators $S_1|H_0$ and $N=N_k|H_0$
satisfy the hypothesis of the previous lemma.

\hspace{15 mm} Let ${\cal H}_m$ be the spaces obtained by applying the lemma.
Since $P_{{\cal H}_1}N|{\cal H}_1$ is invertible its spectral 
radius $r$ is nonzero, so because of the similarity we have
$\norm{P_{{\cal H}_m}N|{\cal H}_m} \geq r >0$ for every $m$,
which contradicts the fact that $N$ is compact. This proves the
theorem.
 
\bigskip

\hspace{15 mm} We remark that from the proof it follows that there
is an obstruction to making the n-tuple either left or right invertible.
 As a consequence of the theorem, if $T$ is  Fredholm
of nonzero index and $k_1,k_2,\cdots ,k_n$ are positive
integers then the n-tuple $(T^{k_1},T^{k_2},\cdots ,T^{k_n})$
has the index equal to zero, but  cannot be perturbed with
compact operators to an invertible n-tuple. The next example
will show that the obstruction to making tuples invertible can
be provided by the index of a subtuple. We still don't know if
this works in general.

\hspace{15 mm} Let ${\bf H}^2({\bf D}^2)$ be the Hardy space on
the bidisk, and $T_{z_1}$ and $T_{z_2}$ the two shifts defined by
$T_{z_1}f(z_1,z_2)=z_1f(z_1,z_2)$, $T_{z_2}f(z_1,z_2)=z_2f(z_1,z_2)$
for $f\in {\bf H}^2({\bf D}^2)$.
It is well known  that the pair $(T_{z_1},T_{z_2})$ is Fredholm of
index 1. Therefore the triple $(T_{z_1},T_{z_2},0)$ is Fredholm of index
zero. Let us show that it cannot be perturbed with compact operators
to a commuting invertible  triple.

\hspace{15 mm} Suppose that there exist compact operators $K_1$, 
$K_2$ and $K_3$ such that the triple $(T_{z_1}+K_1,T_{z_2}+K_2,K_3)$ is
invertible. Let $S_1=T_{z_1}+K_1$ and $S_2=T_{z_2}+K_2$. By theorem 3.8 in
[4] $ind(S_1^n,S_2)=n\cdot ind(S_1,S_2)=n$, which shows that
$dimH^0(S_1^n,S_2)+dimH^2(S_1^n,S_2)\rightarrow \infty $ for $n\rightarrow 
\infty $. So there is a sequence of positive integers $\{n_k\}_k$
such that either $dimH^0(S_1^{n_k},S_2)\rightarrow \infty$ or $dimH^2(S_1^{n_k}
,S_2)\rightarrow \infty $. Without loss of generality we may assume
that $dimH^0(S_1^{n_k},S_2)\rightarrow \infty $. Since $H^0(S_1^n,S_2)=
kerS_1^n\cap kerS_2$ and $kerS_1^n\subset kerS_1^{n+1}$ we get that
$dim(kerS_1^n\cap kerS_2)\rightarrow \infty $ for $n\rightarrow \infty $.
From the spectral mapping theorem it follows that $(S_1^n,S_2,K_3)$
is invertible for any positive integer $n$ hence 
$kerS_1^n\cap kerS_2\cap kerK_3=0$. Therefore we can apply
the lemma to the space $kerS_2$, and to the operators $S_1|kerS_2$ and
$K_3|kerS_2$. Using the same idea as in the proof of the theorem
we contradict the compactness of $K_3$, which proves the claim. 

\vspace{7 mm}
 
\hspace{15 mm} {\bf REFERENCES}

[1] Curto, R., {\em Applications of several complex variables to multiparameter
spectral theory}, Surveys of Some Recent Results in Operator Theory, v.II,
(J. Conway and B. Morrel, eds.), Pitman Res. Notes in Math. Ser. 192,
Longman Sci. and Techn., London, 1988, pp.25--90;

[2] Curto, R., {\em Problems in multivariable operator theory},
 in Contemporary Math.
120 (1991), 15--17;

[3] Gelca, R., {\em Compact perturbations of 
Fredholm n-tuples}, Proceedings AMS,
{\bf 122}(1994), 195--199.

[4] Putinar, M. {\em Some invariants for semi-Fredholm systems of essentialy
commuting operators } JOT 8(1982), 65--90;

[5] Taylor, J., L., {\em A joint spectrum for several commuting operators}, J. 
Funct. Anal. 6(1970), 172--191;

[6] Taylor, J., L., {\em The analytic functional calculus for several 
commuting operators}, Acta Math. 125(1970), 1--38.
\vspace{5 mm}

{\parbox{100 mm}{Department of Mathematics} 
\parbox{70 mm}{Institute of Mathematics}}

{\parbox{100 mm}{The University of Iowa \ \ \ \ \ \ \ \ \ \ \ \ \ \ \ \ \ \ \ and}
\parbox{70 mm}{of the Romanian Academy}}

{\parbox{100 mm}{Iowa City, IA 52242 USA} 
\parbox{70 mm}{P.O.Box 1--764}}

{\parbox{100 mm}{{\em E-mail: rgelca@math.uiowa.edu}} 
\parbox{70 mm}{70700 Bucharest, Romania}}

\medskip

 AMS(MOS) Subj. Classif. 47A53, 47A55, 47A60, 47B05.

\end{document}